%&amstex          
\input amstex\documentstyle{amsppt}  
\pagewidth{12.5cm}\pageheight{19cm}\magnification\magstep1
\NoBlackBoxes
\topmatter
\title Unipotent almost characters of simple $p$-adic groups, II\endtitle
\author G. Lusztig\endauthor
\address{Department of Mathematics, M.I.T., Cambridge, MA 02139}\endaddress
\dedicatory{Dedicated to Professor E. B. Dynkin for his 90th birthday}\enddedicatory
\thanks{Supported in part by National Science Foundation grant DMS-0758262.}\endthanks
\endtopmatter   
\document
\define\ufA{\un{\fA}}

\define\hori{\frac{\hphantom{aaa}}{\hphantom{aaa}}}

\define\dy{\dot y}

\define\mpb{\medpagebreak}

\define\bay{\bar y}

\define\bR{\bar R}

\define\bVV{\bar\VV}
\define\uVV{\un\VV}
\define\baG{\bar\G}

\define\hR{\hat R}

\define\si{\sim}

\define\sqc{\sqcup}

\define\qua{\quad}

\define\hE{\hat E}

\define\dx{\dot x}

\define\tch{\ti\ch}
\define\tSi{\ti\Si}

\define\bZ{\bar Z}

\define\op{\oplus}

\define\part{\partial}

\define\ra{\rangle}
\define\n{\notin}
\define\iy{\infty}
\define\m{\mapsto}
\define\do{\dots}
\define\la{\langle}

\define\lra{\leftrightarrow}

\define\sub{\subset}    

\define\T{\times}
\define\ti{\tilde}
\define\nl{\newline}
\redefine\i{^{-1}}
\define\fra{\frac}
\define\un{\underline}
\define\ov{\overline}
\define\ot{\otimes}

\define\Ad{\text{\rm Ad}}
\define\Hom{\text{\rm Hom}}

\define\Aut{\text{\rm Aut}}

\define\Ker{\text{\rm Ker}}

\define\tr{\text{\rm tr}}

\define\a{\alpha}
\redefine\b{\beta}
\redefine\c{\chi}
\define\g{\gamma}
\redefine\d{\delta}
\define\e{\epsilon}

\redefine\o{\omega}
\define\p{\pi}
\define\ph{\phi}
\define\ps{\psi}
\define\r{\rho}
\define\s{\sigma}
\redefine\t{\tau}
\define\th{\theta}
\define\k{\kappa}
\redefine\l{\lambda}
\define\z{\zeta}
\define\x{\xi}

\redefine\G{\Gamma}
\redefine\D{\Delta}
\define\Om{\Omega}
\define\Si{\Sigma}
\define\Th{\Theta}
\redefine\L{\Lambda}
\define\Ph{\Phi}

\define\dd{\bold d}

\define\kk{\bold k}

\redefine\tt{\bold t}

\define\CC{\bold C}

\define\GG{\bold G}

\define\NN{\bold N}

\define\RR{\bold R}

\define\VV{\bold V}
\define\WW{\bold W}
\define\ZZ{\bold Z}

\define\cg{\Cal G}
\define\ch{\Cal H}

\define\cj{\Cal J}

\define\cm{\Cal M}

\define\co{\Cal O}
\define\cp{\Cal P}

\define\ct{\Cal T}
\define\cu{\Cal U}
\define\cv{\Cal V}

\define\cz{\Cal Z}

\define\ft{\frak t}

\define\fA{\frak A}
\define\fB{\frak B}

\define\fZ{\frak Z}

\define\tg{\ti g}

\define\tit{\ti t}

\define\ty{\ti y}
\define\tA{\ti A}

\define\tH{\ti H}

\define\tO{\ti O}

\define\emp{\emptyset}
\define\KL{KL}
\define\KML{KmL}
\define\LS{LS}
\define\EEE{L1}
\define\ORA{L2}
\define\ICC{L3}
\define\CLA{L4}
\define\CDG{L5}
\define\UNAC{L6}

\subhead 0.1\endsubhead
For any finite group $\G$, a ``nonabelian Fourier transform matrix'' was introduced in \cite{\EEE}.
This is a square matrix whose rows and columns are indexed by pairs formed by an element of $\G$ and an 
irreducible representation of the centralizer of that element (both defined up to conjugation). As shown in
\cite{\ORA}, this matrix, which is unitary with square $1$, enters (for suitable $\G$) in the character 
formulas for unipotent representations of a finite reductive group.

In this paper we extend the definition of the matrix above to the case where $\G$ is a reductive group over 
$\CC$. We expect that this new matrix (for suitable $\G$) relates the characters of unipotent 
representations of a simple $p$-adic group and the unipotent almost characters of that simple $p$-adic group
\cite{\UNAC}. The new matrix is defined in \S1. The definition depends on some finiteness results 
established in 1.2. Several examples are given in 1.4-1.6 and 1.12. In \S2 a conjectural relation with 
unipotent characters of 
$p$-adic groups is stated. In \S3 we consider an example arising from an odd spin group which provides some 
evidence for the conjecture. Since the group $\G$ which enters in the conjecture is described in the
literature only up to isogeny, we give a more precise description for it (or at least for the derived 
subgroup of its identity component) in the Appendix.

{\it Notation.} If $G$ is an affine algebraic group and $g\in G$ we denote by $g_s$ (resp. $g_u$) the 
semisimple (resp. unipotent) part of $g$. Let $\cz_G$ be the centre of $G$ and let $G^0$ be the identity 
component of $G$. Let $G_{der}$ be the derived subgroup of $G$. If $g\in G$ and $G'$ is a subgroup of $G$ we
set $Z_{G'}(g)=\{x\in G';xg=gx\}$.

\head 1. A pairing\endhead
\subhead 1.1\endsubhead
Let $H$ be a reductive (but not necessarily connected) group over $\CC$; let $\Si$ be the set of semisimple 
elements in $H$. Let $[H]$ be the set of irreducible components of $H$. For $x,y$ in $\Si$ we set 
$A_{x,y}=\{z\in H;zxz\i y=yzxz\i\}$; we have $A_{x,y}=\sqc_{h\in[H]}A^h_{x,y}$ where 
$A^h_{x,y}=A_{x,y}\cap h$. Now $Z^0_H(x)\T Z^0_H(y)$ acts on $A_{x,y}$ by $(v,v'):z\m v'zv\i$, leaving
stable each of the subsets $A_{x,y}^h$ $(h\in[H])$.

\proclaim{Lemma 1.2}The $Z^0_H(x)\T Z^0_H(y)$ action on $A_{x,y}$ in 1.1 has only finitely many orbits.
\endproclaim
For any $c\in Z_H(y)\cap\Si$ we choose a maximal torus $T_c$ of $Z_{Z_H(y)}(c)$ and a maximal torus $T'_c$ 
of $Z_H(c)$ that contains $T_c$. (Note that $T_c\sub Z_H(c)$.) We show:

(a) {\it If $E$ is any semisimple $H$-conjugacy class in $H$ then $E\cap cT_c$ is finite.}
\nl
Since $E$ is a finite union of $H^0$-conjugacy classes, it is enough to show that, if $E'$ is any semisimple
$H^0$-conjugacy class in $H$, then $E'\cap cT_c$ is finite. By \cite{\CDG, I, 1.14}, $E'\cap cT'_c$ is 
finite. Since $T_c\sub T'_c$, we see that $E'\cap cT_c$ is also finite, as required.

Let $\{c_i;i\in[1,n]\}$ be a collection of elements of $Z_H(y)\cap\Si$, one in each connected component of 
$Z_H(y)$. Let $\{c'_k;k\in[1,n']\}$ be a collection of elements of $Z_H(x)\cap\Si$, one in each connected 
component of $Z_H(x)$. Let $E$ be the $H$-conjugacy class of $x$. By (a), $E_i:=E\cap c_iT_{c_i}$ is a 
finite set for $i\in[1,n]$. For any $i\in[1,n],e\in E_i$ we set $P_{i,e}=\{g\in H;gxg\i=e\}$. We have 
$P_{i,e}=p_{i,e}Z_H(x)$ for some $p_{i,e}\in H$ hence
$$P_{i,e}=\cup_{k\in[1,n']}p_{i,e}c'_kZ^0_H(x).\tag b$$
Let $z\in A_{x,y}$. We have $zxz\i\in Z_H(y)\cap\Si$ hence $zxz\i\in c_iZ_H^0(y)$ for a unique $i\in[1,n]$.
Using \cite{\CDG, I, 1.14(c)} (with $G$ replaced by the reductive group $Z_H(y)$) we see that $zxz\i$ is 
$Z^0_H(y)$-conjugate to an element of $c_iT_{c_i}$. Thus, $v'zxz\i v'{}\i\in c_iT_{c_i}$ for some 
$v'\in Z^0_H(y)$. Hence $v'zxz\i v'{}\i=e\in E_i$ so that $v'z\in P_{i,e}$. Using (b) we see that 
$v'z=p_{i,e}c'_kv$ for some $v\in Z^0_H(x)$ and some $k\in[1,n']$. Hence $v'zv\i=p_{i,e}c'_k$. We see that 
the finitely many elements $p_{i,e}c'_k$ ($i\in[1,n],k\in[1,n'],e\in E_i$) represent all the orbits in the 
lemma. This completes the proof.

{\it Remark.} The following result can be deduced from the lemma above.

(c) {\it Let $C,C'$ be two semisimple $H^0$-conjugacy classes in $H$. Let $X=\{(g,g')\in C\T C';gg'=g'g\}$. 
Then the $H^0$-action on $X$, $g_1:(g,g')\m(g_1gg_1\i,g_1g'g_1\i)$ has only finitely many orbits.}
\nl
Let $x\in C,y\in C'$. Let $\fZ$ be a set of representatives for the orbits of $Z_{H^0}(x)\T Z_{H^0}(y)$ on 
$A^{H^0}_{x,y}$. This is a finite set. (Even the orbits of the smaller group $Z_H^0(x)\T Z_H^0(y)$ form a 
finite set.) Let $\co$ be an $H^0$-orbit in $X$. In $\co$ we can find an element of the form $(zxz\i,y)$ 
where $z\in H^0$, $zxz\i y=yzxz\i$ hence $z\in A^{H^0}_{x,y}$. Thus $z=v'\z v\i$ where $v'\in Z_{H^0}(y)$, 
$v\in Z_{H^0}(x)$, $\z\in\fZ$. We have $(zxz\i,y)=(v'\z v\i x v\z\i v'{}\i,y)=(v'\z x\z\i v'{}\i,y)\in\co$ 
hence $(\z x\z\i,y)\in\co$. Thus the number of $H^0$-orbits in $X$ is $\le|\fZ|$, proving (c).

\subhead 1.3\endsubhead 
In the remainder of this section we fix a (not necessarily connected) reductive group $H$ over $\CC$ and a 
finite subgroup $\L$ of $\cz_H$. Let $\Si$ be the set of semisimple elements in $H$. 

A pair $x',y'$ of commuting elements in $\Si$ is said to be {\it adapted} if there exists a maximal torus 
$\ct$ of $H^0$ such that $x'\ct x'{}\i=\ct$, $y'\ct y'{}\i=\ct$, $(\ct\cap Z_H(x'))^0$ is a maximal torus 
of $Z_H(x')^0$ and $(\ct\cap Z_H(y'))^0$ is a maximal torus of $Z_H(y')^0$. For example, if $x',y'$ are 
contained in $H^0$ then $x',y'$ is adapted if and only there exist a maximal torus of $H^0$ that contains 
$x'$ and $y'$; this condition is automatically satisfied if $(H^0)_{der}$ is simply connected.

As in 1.1, let $[H]$ be the set of connected components of $H$. Let $x,y\in\Si$ and $h\in[H]$. 

For any $x\in\Si$ we set $\bZ(x)=Z_H(x)/(Z_H^0(x)\L)$. For $x\in\Si$ let $I_x$ be the set of isomorphism 
classes of irreducible representations $\s$ (over $\CC$) of $Z_H(x)$ on which the subgroup $Z_H(x)^0$ acts 
trivially. Let $\hat\L=\Hom(\L,\CC^*)$. If $x\in\Si$ and $\c\in\hat\l$ we denote by $I_x^\c$ the set of 
$\s\in I_x$ such that $\L$ acts on $\s$ through $\c$ times identity. Let $\tSi$ be the set of pairs 
$(x,\s)$ where $x\in\Si$ and $\s\in I_x$. Now $H$ acts on $\tSi$ by $f:(x,\s)\m{}^f(x,\s):=(fxf\i,{}^f\s)$ 
where ${}^f\s\in I_{fxf\i}$ is obtained from $\s$ via $\Ad(f):Z_H(x)@>\si>>Z_H(fxf\i)$.

For any $\c\in\hat\L$ let $\tSi^\c=\{(x,\s)\in\tSi;\s\in I_x^\c\}$. 

We assume given a function $\k:\Si\T\Si\T[H]@>>>\CC$ such that

$\k(fxf\i,f'yf'{}\i,f'hf\i)=\k(x,y,h)$ 

$\k(x,y,h)=\k(y,x,h\i)$ for any 

$\k(\z x,\z'y,h)=\k(x,y,h)$ 
\nl
for any $(x,y,h)\in\Si\T\Si\T[H]$, $f,f'\in H$, $\z,\z'\in\L$.

We say that $\k$ is a {\it weight function.} 
For $(x,y,h)\in\Si\T\Si\T[H]$ let ${}^0A^h_{x,y}$ be a set of representatives for the orbits of the action 
$(v,v'):z\m v'zv\i$ of $Z^0_H(x)\T Z^0_H(y)$ on $A^h_{x,y}$; this is a finite set, by 1.2. 
Let ${}^0\tA^h_{x,y}$ be the set of all $z\in{}^0A^h_{x,y}$ such that $zxz\i,y$ form an adapted pair for $H$.

We say that $\k$ is the standard weight function if $\k(x,y,h)=|{}^0\tA^h_{x,y}|\i$ whenever 
$\tA^h_{x,y}\ne\emp$ and $\k(x,y,h)=0$ whenever $\tA^h_{x,y}=\emp$. This weight function is clearly
independent of the choice of ${}^0A^h_{x,y}$.

We return to the general case. Let $(x,\s)\in\tSi,(y,\t)\in\tSi$. We set
$$\align&((x,\s),(y,\t))=|\L/(\L\cap H^0)|\i|\bZ(x)|\i|\bZ(y)|\i\\&
\T\sum_{h\in[H]}\k(x,y,h)\sum_{z\in{}^0A^h_{x,y}}\ov{\tr(zxz\i,\t)}\tr(z\i yz,\s).\tag a\endalign$$
(Clearly, this is independent of the choice of ${}^0A^h_{x,y}$.)

In the case where $H$ is finite, $\L=\{1\}$ and $\k$ is the standard weight function (which in this case 
satisfies $\k(x,y,h)=1$ whenever $A^h_{x,y}\ne0$)  this reduces to the pairing introduced in \cite{\EEE}. 

Note that we can take ${}^0A^h_{y,x}$ to be the image of ${}^0A^{h\i}_{x,y}$ under $z\m z\i$;
it follows that $((y,\t),(x,\s))=\ov{((x,\s),(y,\t))}$. If 
$(x,\s)\in\tSi^\c,(y,\t)\in\tSi^{\c'}$ with $\c,\c'\in\hat\L$ and if $\z,\z'\in\L$ then
$$((\z x,\s),(\z'y,\t))=\c(\z')\ov{\c'(\z)}((x,\s),(y,\t)).$$
(We can take ${}^0A^h_{\z x,\z'y}={}^0A^h_{x,y}$. For $f,f'\in H$ we have
$$({}^f(x,\s),{}^{f'}(y,\t))=((x,\s),(y,\t)).$$
(We can take ${}^0A^h_{fxf\i,f'yf'{}\i}=f'({}^0A^{f'{}\i hf}_{x,y})f\i$).

Now let $\VV$ be the $\CC$-vector space with basis consisting of the elements $(x,\s)\in\tSi$. Then (a) 
extends to a form $(,):\VV\T\VV@>>>\CC$ which is linear in the first variable, antilinear in the second 
variable and satisfies $(v',v)=\ov{(v,v')}$ for $v,v'\in\VV$. We have $\VV=\op_{\c\in\hat\L}\VV^\c$ where 
$\VV^\c$ is the subspace of $\VV$ spanned by $(x,\s)\in\tSi^\c$. For any $\c,\c'$ in $\hat\L$ let 
$$\cj^\c_{\c'}=\{v\in\VV^\c;(v,\VV^{\c'})=0\},\qua \bVV^\c_{\c'}=\VV^\c/\cj^\c_{\c'}.$$
Then $(,)$ induces a pairing $\bVV^\c_{\c'}\T\bVV^{\c'}_\c @>>>\CC$ (denoted again by $(,)$) which is linear 
in the first variable, antilinear in the second variable and satisfies $(v'_1,v_1)=\ov{(v_1,v'_1)}$ for 
$v'_1\in\bVV^\c_{\c'}$, $v_1\in\bVV^{\c'}_\c$. Let $\fA^\c_{\c'}$ be the image of the $\tSi^\c$ under the 
obvious map $\VV^\c@>>>\bVV^\c_{\c'}$. When $\c=1$ we denote by $(x,\s)_{\c'}$ the image of 
$(x,\s)\in\tSi^1$ under $\VV^1@>>>\bVV^1_{\c'}$.

\subhead 1.4\endsubhead
In 1.4-1.6 we assume that $\k$ is the standard weight function.
In this subsection we assume that $H=H^0\L$ and that $(H^0)_{der}$ is simply connected. In this case any
pair of commuting semisimple elements in $H$ is adapted. Let 
$x,y\in\Si$. We have $\bZ(x)=\{1\}$, $\bZ(y)=\{1\}$, $I_x^1=\{1\}$, $I_y^1=\{1\}$, 
$|[H]|=|\L/(\L\cap H^0)|$, $A^h_{x,y}\ne\emp$ for any $h\in[H]$. From the definitions we see that 
$((x,1),(y,1))=1$. Hence $(x,1)-(x',1)\in\cj^1_1$ for any $x,x'\in\Si$. We see that $\fA^1_1$ consists of a 
single element $(1,1)_1$ which has inner product $1$ with itself and $\dim\bVV^1_1=1$. (This example applies
to the situation in 1.7(a) with $u$ any unipotent element in $\cg=SL_N(\CC)$.)

\subhead 1.5\endsubhead
In this subsection we assume that $H=H^0\sqc H^1$ and that $H^0$ is isomorphic to $\CC^*$ (we denote by 
$\l\m g_\l$ an isomorphism $\CC^*@>\si>>H^0$); we also assume that any $r\in H^1$ satisfies 
$rg_\l r\i=g_{\l\i}$ for all $\l\in\CC^*$ and that $\L=\cz_H$ that is, $\L=\{1,g_{-1}\}\sub H^0$. 
In this case any pair of commuting semisimple elements in $H$ is adapted. We fix 
$r\in H^1$. We have $r^2=1$ or $r^2=g_{-1}$. The case where $r^2=1$ (resp. $r^2=g_{-1}$) arises in the 
situation in 1.7(a) with $u$ a subregular unipotent element in $\cg=Spin_{2n+1}(\CC)$ with $n$ even (resp. 
$n$ odd). We have 
$$Z_H(1)=Z_H(g_{-1})=H, Z_H(g_\l)=H^0\text{ if }\l\in\CC^*-\{1,-1\},$$ 
$$Z_H(r)=\{1,g_{-1},r,rg_{-1}\}.$$
Hence 
$$|\bZ(1)|=2, |\bZ(g_\l)|=1\text{ if }\l\in\CC^*-\{1,-1\}, |\bZ(r)|=2.$$
We have 
$$A_{g_\l,g_{\l'}}=H\text{ for any }\l,\l'\in\CC^*,\qua A_{1,r}=H,\qua A_{g_\l,r}=\emp\text{ for any }
\l\in\CC^*-\{1,-1\},$$
$$A_{r,r}=\{1,g_{-1},r,rg_{-1},g_i,g_{-i},rg_i,rg_{-i}\},\text{ where }i=\sqrt{-1}.$$
Hence we can take
$${}^0A^{H^0}_{g_\l,g_{\l'}}=\{1\},\qua {}^0A^{H^1}_{g_\l,g_{\l'}}=\{r\}\text{ for any }\l,\l'\in\CC^*,$$
$${}^0A^{H^0}_{1,r}=\{1\},\qua {}^0A^{H^1}_{1,r}=\{r\},$$
$${}^0A^{H^0}_{g_\l,r}=\emp,\qua {}^0A^{H^1}_{\l,r}=\emp\text{ for any }\l\in\CC^*-\{1,-1\},$$
$${}^0A^{H^0}_{r,r}=\{1,g_{-1},g_i,g_{-i}\},\qua {}^0A^{H^1}_{r,r}=\{r,rg_{-1},rg_i,rg_{-i}\}.$$
Let $\e$ be the non-trivial $1$-dimensional representation of $H$ which is trivial on $H^0$. The restriction
of $\e$ to $Z_H(r)$ is denoted again by $\e$. From the previous results we can write the $5\T5$ matrix of 
inner products $((x,\s),(y,\t))$ with $(x,\s)$ running through $(1,1),(1,\e),(r,1),(r,\e),(g_\l,1)$ (they 
index the rows), $(y,\t)$ running through $(1,1),(1,\e),(r,1),(r,\e),(g_{\l'},1)$ (they index the columns)
and with $\l,\l'\in\CC^*-\{1,-1\}$:
$$\left(\matrix\fra{1}{2}&\fra{1}{2}&\fra{1}{2}&\fra{1}{2}&1\\
               \fra{1}{2}&\fra{1}{2}&-\fra{1}{2}&-\fra{1}{2}&1\\
               \fra{1}{2}&-\fra{1}{2}&\fra{1}{2}&-\fra{1}{2}&0\\
               \fra{1}{2}&-\fra{1}{2}&-\fra{1}{2}&\fra{1}{2}&0\\
               1         &    1      &    0      &    0     &2\endmatrix\right)$$
From the results in 1.3 we see that $(1,1)_1,(1,\e)_1,(r,1)_1,(r,\e)_1,(g_\l,1)_1$ generate $\bVV^1_1$ and 
from the matrix above we see that for any $\l\in\CC^*-\{1,-1\}$ we have $(g_\l,1)_1=(1,1)_1+(1,\e)_1$. Hence
$(1,1)_1,(1,\e)_1,(r,1)_1,(r,\e)_1$ generate $\bVV^1_1$ and the matrix of their inner products is
$$\left(\matrix\fra{1}{2}&\fra{1}{2}&\fra{1}{2}&\fra{1}{2}\\
               \fra{1}{2}&\fra{1}{2}&-\fra{1}{2}&-\fra{1}{2}\\
               \fra{1}{2}&-\fra{1}{2}&\fra{1}{2}&-\fra{1}{2}\\
               \fra{1}{2}&-\fra{1}{2}&-\fra{1}{2}&\fra{1}{2}\endmatrix\right)$$
which is nonsingular. We see that $\fA^1_1$ consists of $(1,1)_1,(1,\e)_1,(r,1)_1,(r,\e)_1$,
$(1,1)_1+(1,\e)_1$, and the first four of these elements form a basis of $\bVV^1_1$.

\subhead 1.6\endsubhead
Let $E$ be a vector space of dimension $2n\ge2$ over the field $F_2$ with two elements with a given
nondegenerate symplectic form $\la,\ra:E\T E@>>>F_2$. In this subsection we assume that $H$ is a 
Heisenberg group attached to $\la,\ra$ that is, a finite group with a given surjective homomorphism 
$\ps:H@>>>E$ whose kernel consists of $1$ and another central element $c$ with $c^2=1$ and in which for any 
$x,y\in E$ we have $\dx^2=\dy^2=c$, $\dx\dy=c^{\la x,y\ra}\dy\dx$ for any $\dx\in\ps\i(x),\dy\in\ps\i(y)$; 
we also assume that $\L=\{1,c\}$. (This example arises in the situation in 1.7(a) with $u$ a suitable 
unipotent element in $\cg$ of spin type.)

We have $Z_H(1)=Z_H(c)=H$ and for $x\in E-\{0\}$, $\dx\in\ps\i(x)$ we have
$Z_H(x)=\ps\i(\{x'\in E;\la x,x'\ra=0\})$. We denote the nontrivial character of $\L$ by $\c$. If $x\in E$ 
and $\dx\in\ps\i(x)$ then $I_{\dx}^1$ can be identified with $E_x:=E/F_2x$ (to each $\bay\in E_x$ 
corresponds the character $\x\m\la\ps(\x),y\ra$ of $H$ where $y\in E$ represents $\bay$). Now $I_1^\c$ and 
$I_c^\c$ consist of a single representation $\r$ whose character is $2^n$ at $1$, $-2^n$ at $c$ and $0$ at 
all other elements of $H$. If $x\in E-\{0\}$ and $\dx\in\ps\i(x)$ then $I_{\dx}^\c$ consists of two 
representations $\r,\r'$ such that the character of $\r$ at $1,c,\dx,c\dx$ is 
$2^{n-1},-2^{n-1},2^{n-1}i,-2^{n-1}i$ (respectively), the character of $\r'$ at $1,c,\dx,c\dx$ is 
$2^{n-1},-2^{n-1},-2^{n-1}i,2^{n-1}i$ (respectively) and whose character at all other elements of $H$ is 
$0$; here $i=\sqrt{-1}$. Let $x,x',y,y'\in E$, let $\dx\in\ps\i(x)$, $\dx'\in\ps\i(x')$ and let $\bay$ 
(resp. $\bay'$) be the image of $y$ (resp. $y'$) in $E_x$ (resp. $E_{x'}$). From the definitions we have
$$((\dx,\bay),(\dx',\bay'))=0\text{ if }\la x,x'\ra\ne0,$$
$$((\dx,\bay),(\dx',\bay'))=2^{-2n+2}(-1)^{\la x',y\ra+\la x,y'\ra}(\d_{x,0}+1)\i(\d_{x',0}+1)\i$$
if $\la x,x'\ra=0$.

Let $x,x',y\in E$, let $\dx\in\ps\i(x)$, $\dx'\in\ps\i(x')$ and let $\bay$ be the image of $y$ in $E_x$; let 
$\r_1\in I_{\dx'}^\c$. From the definitions we have
$$((\dx,\bay),(\dx',\r_1))=0\text{ if }x\ne0,$$
$$((\dx,\bay),(\dx',\r_1))=2^{-n}(-1)^{\la x',y\ra}\text{ if }\dx=1,$$
$$((\dx,\bay),(\dx',\r_1))=-2^{-n}(-1)^{\la x',y\ra}\text{ if }\dx=c.$$
Let $x,x'\in E$, let $\dx\in\ps\i(x)$, $\dx'\in\ps\i(x')$ and let $\r_1\in I_{\dx}^\c$,
$\r'_1\in I_{\dx'}^\c$. From the definitions we see that $((\dx,\r_1),(\dx',\r'_1))$ is equal $0$ if 
$x\ne x'$, is equal to $1$ if $\dx=\dx',r_1=r'_1$ or if $\dx'=c\dx,r_1\ne r'_1$ and is equal to $-1$ if 
$\dx'=c\dx,r_1=r'_1$ or if $\dx'=\dx,r_1\ne r'_1$.

Let $Z$ be the set of all pairs $(x,\bay)$ where $x\in E$ and $\bay\in E_x$; let $\cm$ be the matrix with 
rows and columns indexed by $Z$ with entries $\cm_{(x,\bay),(x',\bay')}=((\dx,\bay),(\dx',\bay'))$ where 
$\dx\in\ps\i(x)$, $\dx'\in\ps\i(x')$. Its square $\cm^2$ has entries
$$\cm^2_{(x,\bay),(x',\bay')}=\sum_{(x'',\bay'')\in Z}((\dx,\bay),(\dx'',\bay''))((\dx'',\bay''),
(\dx',\bay'))$$
where $\dx\in\ps\i(x)$, $\dx'\in\ps\i(x')$, $\dx''\in\ps\i(x'')$. Let $x,x',y,y'\in E$ and let $\bay$ (resp. 
$\bay'$) be the image of $y$ (resp. $y'$) in $E_x$ (resp. $E_{x'}$). From the definitions we have
$$\align&\cm^2_{(x,\bay),(x',\bay')}=2^{-4n+4}\sum_{(x'',y'')\in E\T E;\la x,x''\ra=0,\la x'',x'\ra=0}\\&
(-1)^{\la x'',y+y'\ra+\la x+x',y''\ra}(\d_{x,0}+1)\i(\d_{x'',0}+1)^{-2}(2-\d_{x'',0})\i(\d_{x',0}+1)\i\\&
=\d_{x,x'}2^{-2n+4}\sum_{x''\in E;\la x,x''\ra=0}(-1)^{\la x'',y+y'\ra}(\d_{x,0}+1)^{-2}(\d_{x'',0}+1)\i 2\i
\\&=\d_{x,x'}2^{-2n+3}\sum_{x''\in E;\la x,x''\ra=0}(-1)^{\la x'',y+y'\ra}(\d_{x,0}+1)^{-2}\\&
-\d_{x,x'}2^{-2n+3}(\d_{x,0}+1)^{-2}+\d_{x,x'}2^{-2n+3}(\d_{x,0}+1)^{-2}2\i\\&
=\d_{x,x'}\d_{\bay,\bay'}2^{-2n+3}|\{x''\in E;\la x,x''\ra=0\}|(\d_{x,0}+1)^{-2}\\&
-\d_{x,x'}2^{-2n+2}(\d_{x,0}+1)^{-2}\\&
=\d_{x,x'}\d_{\bay,\bay'}2^{-2n+3}2^{2n-1}(\d_{x,0}+1)\i-\d_{x,x'}2^{-2n+2}(\d_{x,0}+1)^{-2}\\&
=\d_{x,x'}(2-\d_{x,0})(\d_{\bay,\bay'}2-|E_x|\i).\endalign$$
We see that $\cm^2$ is a direct sum over $x\in E$ of matrices $\cp(x)$ indexed by $E_x\T E_x$ where 
$\cp(x)=(2-\d_{x,0})(2I-|E_x|\i J)$ and all entries of $J$ are $1$. Thus we have $J^2=|E_x|J$. Setting 
$\cp'(x)=(2-\d_{x,0})\i\cp(x)$ we have $|E_x|\i J=2I-\cp'(x)$ hence
$$\cp'(x)^2-4\cp'(x)+4I=(2I-\cp'(x))^2=|E_x|^{-2}J^2=|E_x|\i J=2I-\cp'(x).$$
Thus $\cp'(x)^2-3\cp'(x)+2I=0$ that is $(\cp'(x)-I)(\cp'(x)-2I)=0$. This shows that $\cp'(x)$ is semisimple 
with eigenvalues $1$ and $2$ hence $\cm^2$ is semisimple with eigenvalues $2-\d_{x,0}$ and $2(2-\d_{x,0})$ 
hence $\cm$ is semisimple with eigenvalues $\pm1,\pm\sqrt{2},\pm2$. In particular $\cm$ is invertible.

For $z=(x,\bay)\in Z$ we set $a_z=(\dx,\bay)_1\in\bVV^1_1$ (where $\dx\in\ps\i(x)$); this is independent of 
the choice of $\dx$. Note that $\fA^1_1$ consists of the elements $a_z$. From the fact that $\cm$ is 
invertible we see that $\{a_z;z\in Z\}$ is a basis of $\bVV^1_1$.

\subhead 1.7\endsubhead
Let $\cg$ be a connected, simply connected, almost simple reductive group over $\CC$. Let $u$ be a unipotent
element of $\cg$. Let $\cv$ be the unipotent radical of $Z_\cg(u)^0$.

(a) {\it Until the end of 1.11 we assume that $H=H(u)$ is a reductive subgroup of $Z_\cg(u)$ such that 
$Z_\cg(u)=H\cv$, $H\cap\cv=\{1\}$ and that $\L=\cz_\cg$.}
\nl
(It is well known that such $H$ exists and that $H$ is unique up to $Z_\cg$-conjugacy; note that 
$\cz_\cg\sub H$.) 
We can state the following result in which we allow $\k$ to be any weight function for 
$H$. 

\proclaim{Proposition 1.8} Let $\c,\c'\in\hat\L$. Then $\fA^\c_{\c'}$ is finite. In particular, 
$\dim\bVV^\c_{\c'}<\iy$.
\endproclaim
The proof is given in 1.11 after some preparation in 1.9, 1.10.

\subhead 1.9\endsubhead
For $g\in\cg$ we set $T_g=\cz_{Z_\cg(g_s)}^0$, a torus contained in $Z_\cg(g)$. For $g,g'\in\cg$ we write 
$g\si g'$ if $g'g\i\in T_g=T_{g'}$. This is an equivalence relation on $\cg$. Note that if $g\si g'$ then 
$xgx\i\si xg'x\i$ for any $x\in\cg$. We show:

(a) {\it Let $g,g'\in\cg$ be such that $g\si g'$. Then $Z_\cg(g)=Z_\cg(g')$ and $g_u=g'_u$.}
\nl
The weaker statement that $Z_\cg(g)^0=Z_\cg(g')^0$ is contained in the proof of \cite{\CDG, I, 3.4}. Let 
$s=g_s,u=g_u,s'=g'_s,u'=g_u$. We have $g'=tg$ where $t\in T_g=T_{g'}$. Thus $s'u'=tsu$. Since 
$t\in Z_\cg(g)$ we have $ts=st,tu=ut$. Since $s,t$ are commuting semisimple elements, $ts$ is semisimple and
it commutes with $u$. By uniqueness of the Jordan decomposition of $g'$ we have $s'=ts,u'=u$. Let 
$x\in Z_\cg(g)$. Then $xs=sx$. We have $t\in\cz_{Z_\cg(s)}$. Moreover, $x\in Z_\cg(s)$, hence $xt=tx$ so that
$xs'=xts=tsx=s'x$. Thus $x\in Z_\cg(s')$. We also have $x\in Z_\cg(u)=Z_\cg(u')$. Thus, $x\in Z_\cg(s'u')$ 
that is, $x\in Z_\cg(g')$. Thus, $Z_\cg(g)\sub Z_\cg(g')$. By symmetry we have also $Z_\cg(g')\sub Z_\cg(g)$
hence $Z_\cg(g)=Z_\cg(g')$. This proves (a).

\subhead 1.10\endsubhead
In this subsection we assume that we are in the setup of Proposition 1.8. Let $\p:Z_\cg(u)@>>>H$ be the 
homomorphism such that $\p(x)=x$ for $x\in H$ and $\p(x)=1$ for $x\in\cv$. We show:

(a) {\it Let $(x,\s),(x',\s'),(y,\t)\in\tSi$. Assume that for some $v\in\cv$ we have $vxuv\i\si x'u$. Then 
$Z_H(x)=Z_H(x')$. Assume further that $\s=\s'$ (which makes sense by the previous sentence). Then 
$((x,\s),(y,\t))=((x',\s'),(y,\t))$.}
\nl
From $vxuv\i\si x'u$ and 1.9(a) we see that $Z_\cg(vxuv\i)=Z_\cg(x'u)$. Since $(vxuv\i)_s=vxv\i$,
$(vxuv\i)_u=u$, $(x'u)_s=x'$, $(x'u)_u=u$, we deduce that $Z_{\tH}(vxv\i)=Z_{\tH}(x')$. Using the semidirect
product decompositions $Z_\cg(u)=(vHv\i)\cv$, $Z_\cg(u)=H\cv$ we deduce 
$$Z_{vHv\i}(vxv\i)Z_\cv(vxv\i)=Z_H(x')Z_\cv(x')$$
that is 
$$vZ_H(x)v\i Z_\cv(vxv\i)=Z_H(x')Z_\cv(x').$$
Applying $\p$, we deduce $Z_H(x)=Z_H(x')$. This proves the first assertion of (a). From this we see that for
$z\in H$ we have $z\i yz\in Z_H(x)$ if and only if $z\i yz\in Z_H(x')$. In other words, we have 
$A_{x,y}=A_{x',y}$. It follows that for any connected component $h$ of $H$ we have $A^h_{x,y}=A^h_{x',y}$. 
Let $z\in A_{x,y}=A_{x',y}$. Clearly, $\tr(z\i yz,\s)=\tr(z\i yz,\s')$ (recall that $\s=\s'$). We claim that 
$\tr(zxz\i,\t)=\tr(zx'z\i,\t)$. We set $e:=x'x\i\in H$. From $vxuv\i\si x'u$ we have 
$$\ti e:=x'uvu\i x\i v\i\in T_{x'u}.$$
Hence 
$$z\ti ez\i=ze(xvx\i)v\i z\i\in zT_{x'u}z\i=T'$$
where 
$$T':=z\cz_{Z_\cg(x')}^0z\i=\cz_{Z_\cg(zx'z\i)}^0.$$
Now $yu\in Z_\cg(zx'z\i)$ hence $T'\sub Z_\cg(yu)$ that is, $T'\sub Z_{\tH}(y)=Z_H(y)Z_\cv(y)$. It follows 
that $\p(T')$ is a torus contained in $Z_H(y)$. Hence $\p(T')\sub Z_H^0(y)$. Applying $\p$ to 
$ze(xvx\i)v\i z\i\in T'$ we obtain $zez\i\in\p(T')$ hence $zez\i\in Z_H^0(y)$. We have 
$zx'z\i=(zez\i)(zxz\i)$ where $zez\i\in Z_H^0(y)$ acts trivially on $\t$. Hence
$\tr(zxz\i,\t)=\tr(zx'z\i,\t)$ as claimed. It follows that (a) holds.

Next we show:

(b) {\it Let $\c\in\hat\L$. Let $(x,\s)\in\tSi^\c$, $x'\in\Si$. Assume that for some $g\in\cg$ we have 
$gxug\i\si x'u$. Then there exists $\s'\in I_{x'}^\c$ such that for any $(y,\t)\in\tSi$ we have 
$((x,\s),(y,\t))=((x',\s'),(y,\t))$.}
\nl
From 1.9(a) we have $(gxug\i)_u=u$ hence $gug\i=u$ that is $g\in Z_\cg(u)$. We write $g=vg_1$ with 
$g_1\in H$, $v\in\cv$. We have $g_1xug_1\i=x''u$  where $x''\in\Si$ and $vx''uv\i\si x'u$. From 1.3 we have
$((x'',\s'),(y,\t))=(x,\s),(y,\t))$ where $\s'={}^{g_1}\s\in I_{x''}^\c$. From $vx''uv\i\si x'u$ and (a) we 
have $I_{x''}^\c=I_{x'}^\c$ and $((x'',\s'),(y,\t))=((x',\s'),(y,\t))$. Thus 
$((x,\s),(y,\t))=((x',\s'),(y,\t))$. This proves (b).

\subhead 1.11\endsubhead
We now prove Proposition 1.8. For $g,g'\in\cg$ we write $g\asymp g'$ if $\g g\g\i\si g'$ for some 
$\g\in\cg$. This is an equivalence relation on $\cg$; the equivalence classes are called the strata of 
$\cg$. According to \cite{\CDG, I, 3.7}, $\cg$ has only finitely many strata. Let $S_1,S_2,\do,S_n$ be the 
strata of $\cg$ which have nonempty intersection with $u\Si$. For each $i\in[1,n]$ we choose $x_i\in\Si$ 
such that $ux_i\in S_i$. Now let $(x,\s)\in\tSi^\c$. Let $S$ be the stratum of $\cg$ that contains $xu$. We 
have $S=S_i$ for some $i\in[1,n]$. Hence $gxug\i\si x_iu$ for some $g\in\cg$. By 1.10(b) there exists 
$\s'\in I_{x_i}^\c$ such that $(x,\s)=(x_i,\s')$ in $\bVV^\c_{\c'}$. We see that $\fA^\c_{\c'}$ is the image 
under $V^\c@>>>\bVV^\c_{\c'}$ of the finite set consisting of $(x_i,\s')$ with $i\in[1,n]$, 
$\s'\in I_{x_i}^\c$. The proposition is proved.

\subhead 1.12\endsubhead
In this subsection we assume that $H=PGL_2(\CC)$, $\L=\{1\}$ and that $\k$ is the standard weight function. 
(This example is as in 1.7(a) for a suitable $u$ with $\cg$ of type $F_4$, see A.11.) 

Let $\l\m g_\l$ be an isomorphism of $\CC^*$ onto a maximal torus $T$ of $H$ whose normalizer in $H$ is 
denoted by $N(T)$. Then any semisimple element of $H$ is conjugate to an element of the form $g_\l$, for 
some $\l\in\CC^*$. We have $Z_H(1)=H$, $Z_H(g_{-1})=N(T)$, $Z_H(g_\l)=T$ if $\l\in\CC^*-\{1,-1\}$. Hence 
$|\bZ(1)|=1$, $|\bZ(g_{-1})|=2$, $|\bZ(g_\l)|=1$ if $\l\in\CC^*-\{1,-1\}$. Let $r\in N(T)-T$. We can find 
$\x\in H$ such that $\x g_{-1}\x\i=r$. We have $A_{1,g_\l}=H$ if $\l\in\CC^*$, 
$A_{g_{-1},g_{-1}}=T\cup rT\cup T\x T\cup T\x rT$, $A_{g_{-1},g_\l}=T\cup rT$ if $\l\in\CC^*-\{1,-1\}$, 
$A_{g_\l,g_{\l'}}=T\cup rT$ if '$\l,\l'\in\CC^*-\{1,-1\}$. Hence we can take
$${}^0A^H_{1,g_\l}=\{1\}, {}^0A^H_{g_{-1},g_{-1}}=\{1,r,\x,\x r\},$$
$${}^0A^H_{g_{-1},g_\l}=\{1,r\}\text{ if }\l\in\CC^*-\{1,-1\},$$
$${}^0A^H_{g_\l,g_{\l'}}=\{1,r\}\text{ if }\l\in\CC^*-\{1,-1\}.$$

In our case any pair of commuting semisimple elements is adapted except for a pair $H$-conjugate to the pair
$g_{-1},r$. We have:
$$\k(1,g_\l,H)=1, \k(g_{-1},g_{-1})=1/2,$$
$$\k(g_{-1},g_\l)=1/2\text{ if }\l\in\CC^*-\{1,-1\},$$
$$\k(g_\l,g_{\l'})=1/2\text{ if }\l\in\CC^*-\{1,-1\}.$$
Let $\e$ be the non-trivial $1$-dimensional representation of $N(T)$ which is trivial on $T$. From the 
previous results we can write the $4\T4$ matrix of inner products $((x,\s),(y,\t))$ with $(x,\s)$ running 
through $(1,1),(g_{-1},1),(g_{-1},\e),(g_\l,1)$  (they index the rows), $(y,\t)$ running through 
$(1,1),(g_{-1},1),(g_{-1},\e),(g_{\l'},1)$ (they index the columns) and with $\l,\l'\in\CC^*-\{1,-1\}$:
$$\left(\matrix 1&\fra{1}{2}&\fra{1}{2}&1\\
           \fra{1}{2} &\fra{1}{2}   &  0       &\fra{1}{2}\\
           \fra{1}{2}& 0   &  \fra{1}{2}   &\fra{1}{2}\\
            1       & \fra{1}{2}   & \fra{1}{2}  & 1\endmatrix\right)$$
We see that for any $\l\in\CC^*-\{1,-1\}$ we have $(g_\l,1)=(1,1)$ in $\bVV^1_1$. Also,  
$(1,1)=(g_{-1},1)+(g_{-1},\e)$ in $\bVV^1_1$. Hence $(g_{-1},1),(g_{-1},\e)$ generate $\bVV^1_1$ and the 
matrix of their inner products is
$$\left(\matrix   \fra{1}{2} &  0       \\
                   0   &  \fra{1}{2}\endmatrix \right)$$
which is nonsingular, so that $\fA^1_1$ is equal to $\{(g_{-1},1),(g_{-1},\e)\}$ and is a basis of 
$\bVV^1_1$. 

\head 2. A conjecture\endhead
\subhead 2.1\endsubhead
Let $\kk$ be an algebraic closure of a finite field $F_q$. Let $K=\kk((t))$ where $t$ is an indeterminate. 
Let $K_0=F_q((t))$, a subfield of $K$. Let $G$ be a simple adjoint algebraic group defined and split over 
$F_q$. Let $\cg$ be a connected, simply connected, almost simple reductive group over $\CC$ of type dual to 
that of $G$. Let $P_\emp$ be an Iwahori subgroup of $\GG:=G(K)$ such that $P_\emp\cap G(K_0)$ is an Iwahori 
subgroup of $G(K_0)$. Let $\WW$ be an indexing set for the set of 
$(P_\emp,P_\emp)$-double cosets in $\GG$. We regard $\WW$ as a group as in \cite{\CLA}. As in 
{\it loc.cit.}, $\WW$ is the semidirect product of a (normal) subgroup $\WW'$ (an affine Weyl group) and a 
finite abelian subgroup $\Om$. For any $\x\in\Om$ let ${}^\x\GG$ be the union of the $(P_\emp,P_\emp)$ 
double cosets in $\GG$ indexed by elements in $\WW'\x$. We have $\GG=\sqc_{\x\in\Om}{}^\x\GG$. Let 
$\GG_{rsc}$ be the set of regular semisimple compact elements of $\GG$, see \cite{\UNAC, 3.1}. We have 
$\GG_{rsc}=\sqc_{\x\in\Om}{}^\x\GG_{rsc}$ where ${}^\x\GG_{rsc}={}^\x\GG\cap\GG_{rsc}$. We set 
${}^\x G(K_0)_{rsc}={}^\x\GG_{rsc}\cap G(K_0)$. Let $\cu$ be the set of isomorphism classes of unipotent 
representations of $G(K_0)$, see \cite{\CLA}. 
To each $\r\in\cu$ we attach a sign $\D(\r)\in\{1,-1\}$ as follows: from the definition, to $\r$ corresponds
a unipotent cuspidal representation $\r_0$ of a finite reductive group and $\D(\r)$ is the invariant
$\D(\r_0)$ defined in \cite{\ORA, 4.23, 4.21}. (We almost always have $\D(\r)=1$.)
Recall that \cite{\CLA} gives a bijection $\fZ^1\lra\cu$ where
for $\x\in\Om$, $\fZ^\x$ is the set of all pairs $(g,\s)$ (up to $\cg$-conjugation) where $g\in\cg$ and $\s$
is an irreducible representation of $Z_\cg(g)/Z_\cg(g)^0$ on which $\cz_\cg$ acts according to $\x$
(we can identify $\Om=\Hom(\cz_\cg,\CC^*)$). For $\z\in\fZ^1$ let 
$R_\z$ be the corresponding unipotent representation, let $\hR_\z$ be the corresponding standard 
representation of $G(K_0)$ (which has $R_\z$ as a quotient) and let $\ph_\z^\x$ be $\D(R_\z)$ times
 the restriction of the 
character of $\hR_\z$ to ${}^\x G(K_0)_{rsc}$. For any $\z\in\fZ^\x$ let $\tt_\z:{}^\x G(K_0)_{rsc}@>>>\CC$ 
be the unipotent almost character defined as in \cite{\UNAC, 3.9} (we assume that \cite{\UNAC, 3.8(f)} 
holds).

We now fix a unipotent element $u$ of $\cg$ and a reductive subgroup $H$ of 
$Z_\cg(u)$ as in 1.7. Then the definitions in \S1 can be applied to this $H$ and to $\L=\cz_\cg$ (a subgroup
of $H$). In particular 
$\Si,\tSi,\VV,\bVV^\c_{\c'},\fA^\c_{\c'}$ and the pairing 1.3(b) are defined. If $(x,\s)\in\tSi^1$ then 
$\s$ can be viewed as an irreducible representation of $Z_\cg(g)/Z_\cg(g)^0$ on which $\cz_\cg$ acts 
trivially hence the $\cg$-conjugacy class of $(su,\s)$ can be viewed as an element $\z=\z_{x,\s}\in\fZ^1$.
We write $\hR_{x,\s},\ph_{x,\s}^\x,\tt_{x,\s}$ instead of $\hR_\z,\ph_\z^\x,\tt_\z$. (Note that $\tt_\z$ is 
defined up to multiplication by a root of $1$.) Let $\ufA^1_\x$ be the set of (class) functions  
${}^\x G(K_0)_{rsc}@>>>\CC$ of the form $\ph^\x_{x,\s}$ for some $(x,\s)\in\tSi^1$. Let $\uVV_\x(u)$ be the 
vector space of (class) functions ${}^\x G(K_0)_{rsc}@>>>\CC$ spanned by $\ufA^1_\x$.

\proclaim{Conjecture 2.2} We assume that $\k$ is the standard weight function for $H$.

(a) For any $\x\in\Om$ there is a unique isomorphism of vector spaces 
$\Th_\x:\bVV^1_\x@>\si>>\uVV_\x(u)$ such that for any $(x,\s)\in\tSi^1$ we have 
$\Th_\x((x,\s)_\x)=\ph^\x_{x,\s}$. In particular $\Th_\x$ defines a bijection $\ufA^1_\x@>\si>>\fA^1_\x$.

(b) There is a unique subset $\fB^1_1$ of $\fA^1_1$ such that $\fB^1_1$ is a basis of $\bVV^1_1$ and any 
element of $\fA^1_1$ is an $\RR_{\ge0}$-linear combination of elements in $\fB^1_1$. 

(c) Let $(x,\s)\in\tSi^1$. We set $(x,\s)_1^*=\sum_{b\in\fB^1_1}(b,(x,\s)_1)b\in\bVV^1_1$. We have
$\Th_1((x,\s)^*_1)=\tt_{x,\s}$ up to a nonzero scalar factor.
\endproclaim

\head 3. An example\endhead
\subhead 3.1\endsubhead
We preserve the setup of 2.1 and we assume that $G$ is of type $C_n$ with $n\ge2$; then we can take
$\cg=\text{Spin}_{2n+1}(\CC)$. We take $u$ in 2.1 to be a subregular unipotent element in $\cg$ so that $H$
in 2.1 is as in 1.5 (with $r^2=g_{-1}^n$). We take $\k$ to be the standard weight function for $H$.

The simple reflections of the affine Weyl group $\WW'$ in 2.1 are
denoted by $s_0,s_1,\do,s_n$ where $s_is_{i+1}$ has order $4$ if $i=0$ or $i=n-1$ and order $3$ if 
$0<i<n-1$; $s_i,s_j$ commute if $|i-j|\ge2$. For any $i\in[0,n]$ let $\WW'{}^i$ be the subgroup of $\WW'$ 
generated by $\{s_j;j\in[0,n]-\{i\}\}$; this can be viewed as the Weyl group of the reductive quotient $G^i$
of a maximal parahoric subgroup of $\GG$ containing $P_\emp$. Note that $G^i$ is defined over $F_q$. 

Let $\tch_q$ be the extended affine Hecke algebra (over $\CC$) with parameter $q$ corresponding to $\WW$ and
let $\ch_q$ be unextended affine Hecke algebra with parameter $q$ corresponding to $\WW'$ (a subalgebra of 
$\tch_q$); note that $\tch_q$ contains an element $\o$ such that $\o^2=1$, $\tch_q=\ch_q\op\ch_q\o$ and such
that $x\m\o x\o$ is the algebra automorphism of $\ch_q$ induced by the automorphism of $\WW'$ given by 
$s_i\m s_{n-i}$ for $i\in[0,n]$. For any $i\in[0,n]$ let $\ch^i_q$ be the Hecke algebra with parameter $q$ 
corresponding to $\WW'{}^i$, viewed naturally as a subalgebra of $\ch_q$.

We consider the following $\WW'$-graphs (in the sense of \cite{\KL} but with the $\mu$-function allowed to 
take complex values):
$$v_0\hori v_1\hori v_2\hori\do\hori v_{n-1}\hori v_n,\tag a$$
$$v_1\hori v_2\hori\do\hori v_{n-1},\tag b$$
$$v_0\qua v_n,\tag c$$
$$v_0\hori v_1\hori v_2\hori\do\hori v_{n-1}\hori v_n\hori v'_{n-1}\hori v'_{n-2}\hori\do \hori v'_1.
\tag $d_\l$   $$
(In ($d_\l$), $v'_1$ is also joined with $v_0$.)

In each case the vertex $v_j$ or $v'_j$ has associated set of simple reflections $\{s_j\}$. For (a) we have 
$\mu(v_0,v_1)=2$,  $\mu(v_n,v_{n-1})=2$, all other $\mu$ are equal to $1$. For (b) all $\mu$ are equal to 
$1$. For (c) there are no edges hence no $\mu$-function. For ($d_\l$) with $\l\in\CC^*$ we have 
$\mu(v_0,v_1)=\l,\mu(v_1,v_0)=\l\i$, $\mu(v_n,v'_{n-1})=\l,\mu(v'_{n-1},v_n)=\l\i$, all other $\mu$ are 
equal to $1$. For any $i\in[0,n]$, the $\WW'$-graphs (a),(b),(c) give rise to a $\WW'{}^i$-graph with the 
same vertices (except the one marked with $s_i$), same $\mu$-function and same associated set of simple 
reflections (at any non-removed vertex). This $\WW'{}^i$-graph is denoted by ($a^i$),($b^i$),($c^i$)
respectively. Let $n^i_a,n^i_b,n^i_c$ be the number of vertices marked with $s_i$ in (a),(b),(c)
respectively (this number is $1$ or $0$).

For each of the $\WW'$-graphs (a),(b),(c),$(d_\l)$ we consider the $\CC$-vector space $E^a,E^b,E^c,E^{d_\l}$
spanned by the vertices of the graph, viewed as a $\ch_q$-module in the standard way defined by the 
$\WW'$-graph structure. This $\ch_q$-module extends to a $\tch_q$-module in which $\o$ acts by permuting 
the bases elements according to: $v_i\lra v_{n-i}$ (for (a),(b),(c)) and (in case $(d_\l)$), 
$v_i\lra v'_{n-i}$ if $i\in[1,n-1]$, $v_0\lra v_n$. Let $\hE^a,\hE^b,\hE^c,\hE^{d_\l}$ be the admissible 
representations of $G(K_0)$ generated by their subspace of $P_\emp\cap G(K_0)$-invariant vectors such that 
the natural $\tch_q$-module structure on this subspace is $E^a,E^b,E^c,E^{d_\l}$ respectively. Now the
$\tch_q$-modules $E^a,E^b,E^c,E^{d_\l}$ can be specialized to $q=1$ (using the $\WW'$-graph structure) and
yield $\WW$-modules $E^a_1,E^b_1,E^c_1,E^{d_\l}_1$ respectively. 

If $i\in[0,n]$ for each of the $\WW'{}^i$-graphs ($a^i$),($b^i$),($c^i$) we consider the $\CC$-vector space 
$E^{a^i},E^{b^i},E^{c^i}$ spanned by the vertices of the graph, viewed as a $\ch_q^i$-module in the 
standard way defined by the $\WW'{}^i$-graph structure. The restriction of the $\ch_q$-module $E^a,E^b,E^c$ 
to $\ch^i_q$-module is isomorphic to $E^{a^i}\op\CC^{n^i_a}$, $E^{b^i}\op\CC^{n^i_b}$, 
$E^{c^i}\op\CC^{n^i_c}$ respectively. Here $\CC^m$ is $0$ if $m=0$ while if $m=1$ it is the $\ch^i_q$-module
defined by the $\WW'{}^i$ graph with a single vertex without any associated simple reflection.

Let $\hE^{a^i},\hE^{b^i},\hE^{c^i}$ be the representations of $G^i(F_q)$ generated by their subspace of 
vectors invariant under a Borel subgroup such that the natural $\ch^i_q$-module structure on this subspace 
is $E^{a^i},E^{b^i},E^{c^i}$ respectively. Now the $\ch^i_q$-modules $E^{a^i},E^{b^i},E^{c^i}$ can be 
specialized to $q=1$ (using the $\WW'{}^i$-graph structure) and yield $\WW'{}^i$-modules 
$E^{a^i}_1,E^{b^i}_1,E^{c^i}_1$ respectively. 

With notation in 2.1 and 1.5, $\hE^a,\hE^b,\hE^c$ can be identified with the standard representations 
$\hR_{1,1},\hR_{1,\e},\hR_{r,1}$ of $G(K_0)$. On the other hand the standard representation $\hR_{1,g_\nu}$ 
(where $\nu\in\CC^*-\{1,-1\}$) is among the representations $\hE^{d_\l}$. Let 
$\ph^*_{1,1},\ph^*_{1,\e},\ph^*_{r,1},\ph^*_{1,g_\nu}$ be the restriction of the character of 
$\hR_{1,1},\hR_{1,\e},\hR_{r,1},\hR_{1,g_\nu}$ to the set $G(K_0)_{vrc}$ of very regular compact elements 
(see \cite{\KML}) in $G(K_0)$. According to \cite{\KML}, $G(K_0)_{vrc}$ has a canonical partition
$\sqc_\th G(K_0)^\th_{vrc}$ into pieces invariant under conjugation by $G(K_0)$, where $\th$ runs over the 
set of $\WW$-conjugacy classes of elements of finite order in $\WW'$. Any finite dimensional $\WW$-module 
$\r$ gives rise to a class function $e_\r:G(K_0)_{vrc}@>>>\CC$ where $e_\r|_{G(K_0)^\th_{vrc}}$ is the 
constant equal to the value of $\r$ at any element of $\th$. In particular the class functions 
$e_{E^a_1},e_{E^b_1},e_{E^c_1}$ on $G(K_0)_{vrc}$ are well defined; we denote them by 
$e_{1,1},e_{1,\e},e_{r,1}$. We set $e_{r,\e}=0$ (a function on $G(K_0)_{vrc}$).

We have:
$$\ph^*_{1,g_\nu}=\ph^*_{1,1}+\ph^*_{1,\e}.\tag e$$
To prove (e) it is enough to show that $\hE^{d_\l}$ and $\hE^a\op\hE^b$ have the same character on
$G(K_0)_{vrc}$ for any $\l\in\CC^*$. By the results of \cite{\KML}, the character of $\hE^{d_\l}$ on 
$G(K_0)_{vrc}$ depends only on the restrictions of $E^{d_\l}$ to the various $\ch^i_q$ and in particular is 
independent of $\l$. Hence we can assume that $\l=1$. Thus it is enough to show that 
$E^{d_1}\cong E^a\op E^b$ as a $\ch_q$-module. The $\WW'$-graph $(d_1)$ admits an involution $v_i\lra v'_i$ 
($i\in[1,n-1]$), $v_0\m v_0$, $v_n\m v_n$. Its eigenspaces give the required direct sum decomposition of the 
$\ch_q$-module $E^{d_1}$.

We have the following equalities.
$$\ph^*_{1,1}=1/2(e_{1,1}+e_{1,\e}+e_{r,1}+e_{r,\e}),\tag f$$
$$\ph^*_{1,\e}=1/2(e_{1,1}+e_{1,\e}-e_{r,1}-e_{r,\e}),\tag g$$
$$\ph^*_{r,1}=1/2(e_{1,1}-e_{1,\e}+e_{r,1}-e_{r,\e}).\tag g$$
To prove (f) it is enough to show that for any $i\in[0,n]$ and any conjugacy class $\th_0$ of $\WW'{}^i$, the
two sides of (f) take the same value on $G(K_0)^\th_{vrc}$ where $\th$ is the $\WW$-conjugacy class in 
$\WW'$ that contains $\th_0$. By \cite{\KML}, $\ph^*_{1,1}|_{G(K_0)^\th_{vrc}}$ is the constant equal to 
$n^i_a$ plus the value of the character of $\hE^{a^i}$ at a regular semisimple element of $G^i(F_q)$ of type 
$\th_0$. By the results of \cite{\ORA}, this is equal to $n^i_a$ plus the character of
$1/2(E^{a^i}_1+E^{b^i}_1+E^{c^i}_1)$ at $\th_0$ hence it is equal to $n^i_a$ plus the character of 
$1/2(e_{1,1}-n^i_a+e_{1,\e}-n^i_b+e_{r,1}-n^i_c+e_{r,\e})$ at $\th$. It remains to note that
$n^i_a=1/2(n^i_a+n^i_b+n^i_c)$. The proof of (g) and (h) is entirely similar; it is based on \cite{\KML},
\cite{\ORA} and the equalities $n^i_b=1/2(n^i_a+n^i_b-n^i_c)$, $n^i_c=1/2(n^i_a-n^i_b+n^i_c)$.

Now $\hR_{r,\e}$ is an irreducible representation of $G(K_0)$ which is not Iwahori-spherical. Let
$\ph^*_{r,\e}$ be the restriction of its character to $G(K_0)_{vrc}$. Using again \cite{\KML}, \cite{\ORA},
we see that
$$\ph^*_{r,\e}=1/2(e_{1,1}-e_{1,\e}-e_{r,1}+e_{r,\e}).\tag i$$
\nl
From (f)-(i) we deduce:
$$1/2(\ph^*_{1,1}+\ph^*_{1,\e}+\ph^*_{r,1}+\ph^*_{r,\e})=e_{1,1},$$
$$1/2(\ph^*_{1,1}+\ph^*_{1,\e}-\ph^*_{r,1}-\ph^*_{r,\e})=e_{1,\e},$$
$$1/2(\ph^*_{1,1}-\ph^*_{1,\e}+\ph^*_{r,1}-\ph^*_{r,\e})=e_{r,1},$$
$$1/2(\ph^*_{1,1}-\ph^*_{1,\e}-\ph^*_{r,1}+\ph^*_{r,\e})=e_{r,\e}.$$
which, together with (e), confirm the Conjecture 2.2 in a very special case.

\head Appendix \endhead   
\subhead A.0\endsubhead 
Let $\cg,u,\cv,H=H(u)$ be as in 1.7. In this appendix we give some information on the structure of $H$.

\subhead A.1\endsubhead
Assume first that $\cg=SL(V)$ where $V$ is a $\CC$-vector space of dimension $N\ge2$. We assume that 
$V=\op_{i\ge1}V_i$ where $V_i=W_i\ot E_i$ and $W_i,E_i$ are $\CC$-vector spaces of dimension $i,m_i$ 
respectively. Let 
$I=\{i\ge1;m_i\ge1\}$. For $i\in I$ we have an imbedding $\t_i:GL(E_i)@>>>GL(V_i)$, $g\m 1_{W_i}\ot g$. Let 
$\Ph:\prod_{i\in I}GL(E_i)@>>>GL(V)$ be the homomorphism $(g_i)\m\op_{i\in I}\t_i(g_i)$ (an imbedding). Let 
$H=SL(V)\cap\text{ image of }\Ph$. Then $H$ is of the form $H(u)$ for some $u\in\cg$ (as in 1.7) and all
$H(u)$ as in 1.7 are obtained up to conjugacy. Note that $\Ph\i(H)$ is the set of all $(x_i)_{i\in I}$ with 
$x_i\in GL(E_i)$, $\prod_{i\in I}\det(x_i)^i=1$ and $\Ph\i((H^0)_{der})$ is $\prod_{i\in I}SL(E_i)$. Thus
$(H^0)_{der}$ is simply connected.

\subhead A.2\endsubhead
For any $\CC$-vector space $V$ with a fixed symmetric or symplectic nondegenerate bilinear form 
$(,):V\T V@>>>\CC$ we denote by $Is_V$ be the group of isometries of $(,)$.

Until the end of A.8 we assume that $V$ is a $\CC$-vector space of dimension $N\ge1$ with a fixed symmetric
nondegenerate bilinear form $(,):V\T V@>>>\CC$; we write $O_V$ instead of $Is_V$. Let 
$V_*=\{v\in V;(v,v)=1\}$. For any 
$v\in V_*$ define $r_v\in O_V$ by $r_v(v)=-v$, $r_v(v')=v'$ if $v'\in V,(v,v')=0$ (a reflection). 

Let $C(V)$ be the Clifford algebra of $(,)$ that is, the quotient of the tensor algebra of $V$ by the 
two-sided ideal generated by the elements $v\ot v'+v'\ot v-2(v,v')$ with $v,v'\in V$. As a vector space we 
have $C(V)=C(V)^+\op C(V)^-$ where $C(V)^+$ (resp. $C(V)^-$) is spanned by the products $v_1v_2\do v_n$ with
$v_i\in V$ and $n$ even (resp. $n$ odd). Note that for $v\in V_*$ we have $v^2=1$ in $C(V)$. Let $\tO_V$ be 
the subgroup of the group of invertible elements of $C(V)$ consisting of the elements 
$$v_1v_2\do v_n\text{ with }v_i\in V_*, n\ge0.\tag a$$
Let ${}^0\tO_V=\tO(V)\cap C(V)^+$, ${}^1\tO_V=\tO_V\cap C(V)^-=\tO_V-{}^0\tO_V$. Then ${}^0\tO_V$ is a 
subgroup of index $2$ of $\tO_V$. We set $c=-1\in C(V)$. If $v\in V$, $(v,v)=-1$, we have $v^2=-1$ in 
$C(V)$. Hence $c\in{}^0\tO_V$. If $\x\in\tO_V,v\in V$, we have $\x v\x\i\in V$. (Indeed, for 
$\x=v_1v_2\do v_n$ as in (a), we have $\x v\x\i=(-1)^nr_{v_1}r_{v_2}\do r_{v_n}$.) Hence $\x\m[v\m\x v\x\i]$
is a homomorphism $\b:\tO_V@>>>O_V$ with image $O_V$, if $N$ is even, and $O_V^0$ if $N$ is odd. The kernel 
of $\b:{}^0\tO_V@>>>O_V^0$ is $\{1,c\}$. If $N\ge2$ we have ${}^0\tO_V=\tO_V^0$. If $N=1$, ${}^0\tO_V$ is of
order $2$.

\subhead A.3\endsubhead
In the setup of A.2 we assume $\cg={}^0\tO_V$. We also assume (until the end of A.8) that $V=\op_{i\ge1}V_i$
where $V_i=W_i\ot E_i$ and $W_i,E_i$ are $\CC$-vector spaces of dimension $i,m_i$ respectively with given 
nondegenerate bilinear forms $(,)$ (both symmetric if $i$ is odd, both symplectic if $i$ is even) such that 
$(w\ot e,w'\ot e')=(w,w')(e,e')$ for $w,w'\in W_i$, $e,e'\in E_i$ and $(V_i,V_j)=0$ for $i\ne j$. 

Let $I=\{i\ge1;m_i\ge1\}$, $I_{odd}=I\cap(2\ZZ+1)$, $I_{even}=I\cap(2\ZZ)$. For any $t\ge2$ let 
$I^{\ge t}_{odd}=\{i\in I_{odd};m_i\ge t\}$. For any $i\in I$ let 
$\b_i:\tO_{V_i}@>>>O_{V_i}$ be the homomorphism defined like $\b$ with $V$ replaced by $V_i$. For any 
$i\in I_{odd}$ let $w_1^i,w_2^i,\do,w^i_i$ be an orthonormal basis of $W_i$ with respect to $(,)$.

In this appendix, any product over $I$ or $I_{odd}$ is taken using the order of $I$ or $I_{odd}$ induced
from the obvious order on $\NN$. The imbedding $\t_i:Is_{E_i}@>>>O_{V_i}$, $g\m 1_{W_i}\ot g$, restricts to 
an imbedding $Is^0_{E_i}@>>>O^0_{V_i}$. Let $\Ph:\prod_{i\in I}Is_{E_i}@>>>O_V$ be the homomorphism 
$(g_i)\m\op_{i\in I}\t_i(g_i)$ (an imbedding); it restricts to an imbedding 
$\Ph^0:\prod_{i\in I}Is^0_{E_i}@>>>O^0_V$. Let $\bR=O^0_V\cap\text{ image of }\Ph$. We have 
$\text{ image of }\Ph^0=\bR^0$.

For $i\in I_{odd}$ we choose $e_i\in E_{i*}$ and we define $y_i\in O_V$ by $y_i=1$ on $W_i\ot\CC e_i$ and 
$y_i=-1$ on the perpendicular to $W_i\ot\CC e_i$ in $V$. Thus $y_i=\t_i(-r_{e_i})$ on $V_i$, $y_i=-1$ on 
$V_j$ for $j\ne i$. (Here $r_{e_i}$ is a reflection in $E_i$.) We see that $y_i\in\text{ image of }\Ph$,
$y_i^2=1$. 

Let $g=\op_{j\in I}\t_j(\ft_j)\in O_V$ where $\ft_j\in Is_{E_j}$. From the definitions we see that for any 
$i\in I_{odd}$ we have

(a) $y_igy_i\i=\op_{j\in I}\t_j(\ft'_j)\in O_V$ where $\ft'_j\in Is_{E_j}$ is given by
$\ft'_i=r_{e_i}\ft_ir_{e_i}$ and $\ft'_j=\ft_j$ if $j\ne i$. In particular if $g\in\bR^0$ then 
$y_igy_i\i\in\bR^0$.
\nl
Let $\baG$ be the (finite abelian) subgroup of $O_V$ generated by $\{y_i;i\in I_{odd}\}$. Let $\baG^+$ be 
the subgroup of $\baG$ consisting of elements which are products of an even number of generators $y_i$. If 
$I_{odd}=\emp$ we have $\baG=\baG^+=\{1\}$. If $I_{odd}\ne\emp$ then $\baG^+$ is a subgroup of index $2$ of 
$\baG$. We have $\bR=\bR^0\baG^+$ and $\bR^0\cap\baG^+=\{1\}$. Let $H=\{g\in\cg;\b(x)\in\bR\}$. Then $H$ is 
of the form $H(u)$ for some $u\in\cg$ (as in 1.7) and all $H(u)$ as in 1.7 are obtained up to conjugacy. 

\subhead A.4\endsubhead
Our next objective is to describe the structure of $H$, see A.6(e), A.7.

For $i\in I$ we set $S_i={}^0\tO_{E_i}$ if $i\in I_{odd}$ and $S_i=Is_{E_i}$ if $i\in I_{even}$. For 
$i\in I$ let $R'_i=\{\x\in{}^0\tO_{V_i};\b_i(\x)\in\t_i(Is^0_{E_i})\}$. We set $R_i=R'_i$ if $i\in I_{odd}$ 
and $R_i=R'_i{}^0$ if $i\in I_{even}$. 

For $i\in I_{odd}$ let $c'_i$ be $-1$, viewed as an element of $S_i$. For $i\in I$ let $c_i$ be $-1$, viewed
as an element of ${}^0\tO(V_i)$; note that, if $i\in I_{odd}$, we have $c_i\in R_i$.

If $i\in I_{odd}$ then $R_i$ contains
$$x_{i;e,f}=(w_i^1\ot e)\do(w_i^i\ot e)(w_i^1\ot f)\do(w_i^i\ot f)=-x_{i;f,e}$$ 
for any $e,f\in E_{i*}$. (Indeed, $x_{i;e,f}$ projects to $(1_{W_i}\ot r_e)(1_{W_i}\ot r_f)$, where 
$r_e,r_f$ are reflections in $E_i$.) Since $r_er_f$ generate $\t_i(O^0_{E_i})$, we see that the elements 
$x_{i;e,f}$ generate $R_i$. 

If $i\in I_{odd}$ and $m_i=1$ we have $\t_i(Is^0_{E_i})=\{1\}$ and $R_i=\{1,c_i\}$; we see that there is a 
unique isomorphism $S_i@>\si>>R_i$. Thus $\t^0_i$ lifts uniquely to an isomorphism $\tit_i:S_i@>\si>>R_i$ 
which carries $c'_i$ to $c_i$.

If $i\in I_{odd}^{\ge2}$, from an argument in \cite{\ICC, 14.3}, we see that $R_i$ is connected; being a 
double covering of $O^0_{E_i}$, there is a unique isomorphism ${}^0\tO_{E_i}@>\si>>R_i$ which induces the 
identity on $O^0_{E_i}$. Thus $\t^0_i$ lifts uniquely to an isomorphism $\tit_i:S_i@>\si>>R_i$ which carries
$c'_i$ to $c_i$.

If $i\in I_{even}$ then, since $R'_i$ is a double covering (with kernel $\{1,c_i\}$) of $S_i$ which is 
simply connected, we see that $R'_i=R_i\sqc c_iR_i$ and that $\t_i$ lifts uniquely to an isomorphism 
$\tit_i:S_i@>\si>>R_i$.

We set $S=\prod_{i\in I}S_i$; for any $s\in S$ we write $s_i$ for the $i$-component of $s$. Putting together
the isomorphisms $\tit_i$ we get 
$$S@>\si>>\prod_{i\in I}R_i\sub\prod_{i\in I}{}^0\tO_{V_i},$$
Composing with the homomorphism 
$$\prod_{i\in I}{}^0\tO_{V_i}@>>>\cg\tag a$$ 
induced by $\prod_{i\in I}C(V_i)@>>>C(V)$ (multiplication in $C(V)$) we obtain the homomorphism $\Ph'$ in 
the commutative diagram
$$\CD S@>\Ph'>>\cg\\
@V\a VV         @V\b VV\\
\prod_{i\in I}Is^0_{E_i}@>\Ph^0>>O^0_V\endCD\tag b$$
where $\a$ is the obvious homomorphism.

Let $R$ be the image of $\Ph'$. Clearly, $R$ is a closed subgroup of $\cg$. We show that $\b(R)=\bR^0$.
Since $\a$ is surjective, this follows from (b) and the fact that the image of $\Ph^0$ is $\bR^0$. Let 
$$A=\{s\in S;s_i=1\text{ if }i\in I_{even},s_i=c'_i{}^{n_i}\text{ if }i\in I_{odd},
n_i\in\ZZ/2,\sum_{i\in I_{odd}}n_i=0\}.$$
We show:
$$\Ker(\Ph')=A.\tag c$$
From the definitions, if $s\in S$ is given by $s_i=1$ if $i\in I_{even}$, $s_i=c'_i{}^{n_i}$ if 
$i\in I_{odd}$ where $n_i\in\ZZ/2$, then $\Ph'(s)=c^m$ where $m=\sum_{i\in I_{odd}}n_i$. In particular, we 
have $A\sub\Ker(\Ph')$. Conversely, let $(\x_i)_{i\in I}\in S$ be an element in $\Ker(\Ph')$. We set 
$\ti\x_i=\tit_i(\x_i)\in R_i$. Viewing $\ti\x_i$ as an element of $C(V_i)\sub C(V)$, we have 
$\prod_{i\in I}\ti\x_i=1$ (product in $C(V)$). 
Let $I_0$ be the set of all $i\in I$ such that $\ti\x_i$ is a (nonzero) scalar. From the definition of 
$C(V)$ we see that, if $b_i^j$ $(j\in[1,2^{\dim V_i}]$ is a basis of $C(V_i)$, then the products 
$b^f:=\prod_{i\in I}b_i^{f(i)}$ (for various functions $f$ which to each $i\in I$ associate 
$f(i)\in[1,\dim C(V_i)]$) form a basis of $C(V)$. We can assume that $b_i^1=\ti\x_i$ for all $i\in I$ and 
$b_i^2=1$ for all $i\in I-I_0$. Let $f_1,f_2$ be such that $f_1(i)=1$ for all $i\in I$, $f_2(i)=2$ for 
$i\in I-I_0$, $f_2(i)=1$ for $i\in I_0$. We have $b^{f_1}=1$, $b^{f_2}=\l 1$ with $\l\in\CC^*$. Thus 
$b^{f_2}=\l b^{f_1}$. Since $b^{f_2},b^{f_1}$ are part of a basis of $C(V)$, we see that $f_1=f_2$ hence 
$I-I_0=\emp$. Thus $\ti\x_i$ is a nonzero scalar for any $i\in I$. It follows that the element 
$\ti\x_i\in{}^0\tO_{V_i}$ satisfies $\b_i(\ti\x_i)=1$ hence $\ti\x_i=c_i^{n_i}$ (hence $\x_i=c'_i{}^{n_i}$) 
for some $n_i\in\ZZ/2$. If in addition, $i\in I_{even}$, then $\ti\x_i\in R_i$ and $c_i\n R_i$ implies that 
$n_i=0$. Now if $i\in I_{odd}$, then $c_i$ viewed as an element of $C(V)$ is equal to $c$. Hence from 
$\prod_{i\in I}\ti\x_i=1$ (in $C(V)$) it follows that $\sum_{i\in I_{odd}}n_i=0$. Thus $\Ker(\Ph')\sub A$. 
The opposite containment is obvious. This proves (c).

We show:

(d) {\it If $I^{\ge2}_{odd}\ne\emp$ then $R=R^0=H^0$ and $c\in H^0$.}
\nl
If $i\in I$ and $m_i\ge2$, then $R_i=R_i^0$. Hence 
$$\prod_{i\in I}R_i=\prod_{i\in I;m_i=1}R_i\T\prod_{i\in I;m_i\ge2}R_i^0.$$
Thus $R$ is generated by $R^0$ and by $\Ph'(c'_i)$ for various $i\in I$ such that $m_i=1$. To show that
$R=R^0$ it is enough to show that if $i\in I$ is such that $m_i=1$ (hence $i$ is odd) then 
$\Ph'(c'_i)\in R^0$. By assumption we can find $i'\in I_{odd}^{\ge2}$. From (c) we see that 
$\Ph'(c'_ic'_{i'})=1$. Hence 
$$\Ph'(c'_i)=\Ph'(c'_{i'})\in\Ph'(R_{i'})=\Ph'(R_{i'}^0)\sub R^0.$$
Thus we have $R=R^0$. Since $\b(R)=\bR^0$, $H=\b\i(\bR)\cap\cg$, we have $\dim R=\dim H=\dim\bR$ and 
$R\sub H$ so that $R^0=H^0$. From the proof of (c) we see that, if $s\in S$ is given by $s_i=c'_i$ for
$i=i'$ (as above), $s_i=1$ for $i\in I-\{i'\}$, then $\Ph'(s)=c$. Thus, $c\in R$. Since $R=H^0$ we see that
$c\in H^0$; (d) is proved.

We show:

(e) {\it If $I^{\ge2}_{odd}\ne\emp$ then $\b$ induces an isomorphism $H/H^0@>\si>>\bR/\bR^0=baG^+$. Hence 
$H/H^0$ is an finite abelian group of order $2^{|I_{odd}-1|}$.}
\nl
From A.3 we have $\bR/\bR^0=\baG^+$. From the definitions, $\b$ induces a surjective homomorphism 
$H@>>>\bR$ with kernel $\{1,c\}$. It is then enough to note that $\{1,c\}\sub H^0$ (see (d)).

\subhead A.5\endsubhead
Let $i\in I_{odd}$. We define an automorphism $a_i:S@>>>S$ by $s\m s'$ where $s'_i=e_is_ie_i\i$ (product in 
$\tO_{E_i}$), $s'_j=s_j$ for $j\ne i$. Since $a_i$ induces the identity map on $\Ker(\Ph')$ (see A.4(c)), it
follows that $a_i$ induces an automorphism $a'_i:R@>>>R$ (recall that $R$ is the image of $\Ph'$). We set
$$\ty_i=(w_i^1\ot e_i)\do(w_i^i\ot e_i)\in\tO_V.$$
We show:

(a) {\it For any $\tg\in R$ we have $\ty_i\tg\ty_i\i=a'_i(\tg)\in R$.}
\nl
Assume first that $m_i=1$. We have $a_i=1$ (in this case, $e_i$ is in the centre of $O_{E_i}$). Now $R_i$ is
generated by $x_{i;e_i,e_i}=\ty_i^2$ (see A.4) and by $c_i$ hence $\ty_i$ commutes with any element of 
$R_i$; it also commutes with any element of $R_j,j\ne i$ hence $\ty_i\tg\ty_i\i=\tg$. Thus (a) is proved in 
this case. Next we assume that $m_i\ge2$. We have $\b(\ty_i)=y_i$ hence 
$$\b(\ty_i\tg\ty_i\i)=y_i\b(\tg)y_i\i\in y_i\bR^0y_i\i=\bR^0$$
(see (4.3(a))). Thus, $\ty_i R\ty_i\i\in\b\i(\bR^0)\cap\cg\sub H$. Since $R=H^0$, see A.4(d), we see that 
$\ty_i R\ty_i\i=\ty_i H^0\ty_i\i$ is a connected subgroup of $H$. It follows that $\ty_i R\ty_i\i=H^0=R$. 
Thus $Ad(\ty_i)$ is an automorphism of $R$. From the definition and A.3(a) we have
$\b(a'_i(\tg))=y_i\b(\tg)y_i\i=\b(\ty_i\tg\ty_i\i)$. Thus $\ty_i\tg\ty_i\i=a'_i(\tg)f(\tg)$ where 
$f:R@>>>\ker\b$ is a morphism. Now $R$ is connected (see A.4(d)) and $\ker\b$ is finite hence $f$ is 
constant. Clearly, $f(1)=1$. It follows that $f(\tg)=1$ for any $\tg\in R$, proving (a).

\subhead A.6\endsubhead
Let $\D$ be the group defined by the generators $z_i (i\in I_{odd})$ and $c_*$ with relations: $c_*^2=1$, 
$c_*z_i=z_ic_*$, $z_i^2=c_*^{i(i-1)/2}$ for $i\in I_{odd}$ and $z_iz_{i'}=c_*z_{i'}z_i$ for $i\ne i'$ in 
$I_{odd}$. (If $I_{odd}=\emp$, $\D$ is the group of order $2$ with generator $c_*$.) From the definition we 
see that any element of $\D$ can be written in the form $c_*^n\prod_{i\in I_{odd}}z_i^{n_i}$ where 
$n_i\in\ZZ/2$, $n\in\ZZ/2$. It follows that $|\D|\le 2^{|I_{odd}|+1}$. The assignment $c_*\m c$, 
$z_i\m\ty_i=(w_i^1\ot e_i)\do(w_i^i\ot e_i)\in\tO_V$ defines a homomorphism $\z:\D@>>>\tO_V$. For any 
$i\in I_{odd}$ we have $\b(\ty_i)=y_i$. Thus $\b(\z(\D))$ is the subgroup of $O_V$ generated by 
$y_i (i\in I_{odd})$, an (abelian) group of order $2^{|I_{odd}|}$. Since $\z(\D)$ contains $c$ which is in 
the kernel of $\b$, we deduce that $|\z(\D)|\ge2^{|I_{odd}|+1}$ hence $|\D|=|\z(\D)|=2^{|I_{odd}|+1}$. In 
particular, $\z$ defines an isomorphism $\D@>\si>>\z(\D)$ and any element of $\D$ can be written uniquely in
the form $c_*^n\prod_{i\in I_{odd}}z_i^{n_i}$ where $n_i\in\ZZ/2$, $n\in\ZZ/2$.

We now define a homomorphism $\D@>>>\Aut(S)$, $\d\m[s\m{}^\d s]$ by $c_*\m 1$, $z_i\m a_i$ where $a_i$ is as
in A.5. We form the semidirect product $S\cdot\D$ in which, for $s,s'\in S,\d,\d'\in\D$, we have
$(s\d)(s'\d')=(s{}^\d s')(\dd')$. We extend $\z:\D@>>>\tO_V$ to a homomorphism $S\cdot\D@>>>\tO_V$ (denoted 
again by $\z$) by definining $\z$ on $S$ to be $\Ph'$. To show that this is well defined we must verify for 
$i\in I_{odd},s\in S$ that $\Ph'(a_i(s))=\ty_i\Ph'(s)\ty_i\i$. This follows from A.5(a). We show:

(a) {\it The kernel of $\z:S\cdot\D@>>>\tO_V$ is equal to $A'=\{1\}$ if $I_{odd}=\emp$ and to}
$$\align&A':=\{sc_*^n;s\in S,n\in\ZZ/2,s_i=c'_i{}^{n_i}\text{ if }i\in I_{odd},s_i=1\text{ if }i\in I_{even},
\\&n_i\in\ZZ/2,n=\sum_{i\in I_{odd}}n_i\}\endalign$$
{\it if $I_{odd}\ne\emp$.}
\nl
Let $s\in S,\d\in\D$ be such that $s\d\in\Ker(\z)$. We write $\d=c_*^n\prod_{i\in I_{odd}}z_i^{u_i}$ where 
$u_i\in\ZZ/2$, $n\in\ZZ/2$. We set $\x_i=\tit_i(s_i)\in R_i$.

Assume first that $I_{odd}\ne\emp$. We have $1=c^n\prod_{i\in I}\x_i\prod_{i\in I_{odd}}\ty_i^{u_i}$ hence 
$\prod_{i\in I}f_i=\text{power of }c$, where $f_i=\ti\x_i\ty_i^{m_i}$ for $i\in I_{odd}$, $f_i=\x_i$ for 
$i\in I_{even}$. We have $f_i\in C(V_i)^{u_i}$ if $i\in I_{odd}$ and $f_i\in C(V_i)^0$ if $i\in I_{even}$. 
Now the subspaces $\prod_{i\in I}C(V_i)^{r_i}$ of $C(V)$ (with $r_i\in\{0,1\}$) form a direct sum 
decomposition of $C(V)$. It follows that $u_i=0$ for any $i\in I_{odd}$. Thus we have 
$c^n\prod_{i\in I}\x_i=1$. If $n=0$ then from the previous equality we deduce as in the proof of A.5(c) 
that $s_i=c'_i{}^{n_i}$ if $i\in I_{odd}$, $s_i=1$ if $i\in I_{even}$ where $n_i\in\ZZ/2$, 
$\sum_{i\in I_{odd}}n_i=0$. If $n=1$ we pick $i_0\in I_{odd}$ and we define $s'\in S$ by 
$s'_{i_0}=c'_{i_0}s_{i_0}$, $s'_i=s_i$ if $i\in I-\{i_0\}$. Then $\z(s')=1$ from which we deduce as above 
that $s'_i=c'_i{}^{n_i}$ if $i\in I_{odd}$ with $n_i\in\ZZ/2,\sum_{i\in I_{odd}}n_i=0$. Hence 
$s'_i=c'_i{}^{n'_i}$ for $i\in I_{odd}$, $s_i=1$ for $i\in I_{ev}$ with $n'_i\in\ZZ/2$, 
$\sum_{i\in I_{odd}}n'_i=1$. We see that $\ker\z\sub A'$. The opposite containment is obvious. This 
completes the proof when $I_{odd}\ne\emp$.

Assume next that $I_{odd}=\emp$. We have $\D=\{1,c_*\}$. In our case the homomorphism $\a$ in A.4(b) defines
a connected covering of a simply connected group hence is an isomorphism; since $\Ph^0$ in A.4(b) is 
injective it follows that $\b$ in A.4(b) is injective when restricted to $R$, hence $c\n R$. We have 
$1=c^n\prod_{i\in I}\x_i$. Hence $c^n\in R$. Since $c\n R$ it follows that $n=0$. Thus, 
$1=\prod_{i\in I}\x_i$. From this we deduce as in the proof of A.4(c) that $s=1$. This completes the 
proof of (a).

\mpb

Let $\D^+$ be the subgroup of $\D$ consisting of the elements of the form 
$c_*^n\prod_{i\in I_{odd}}z_i^{n_i}$ where $n_i\in\ZZ/2$, $n\in\ZZ/2$, $\sum_{i\in I_{odd}}n_i=0$. We have 
$|\D^+|=2^{|I_{odd}|}$ if $I_{odd}\ne\emp$, $|\D^+|=2$ if $I_{odd}=\emp$. We regard the semidirect product 
$S\cdot\D^+$ as a subgroup of $S\cdot\D$ in an obvious way. Let $\z_0:S\cdot\D^+@>>>\cg$ be the restriction
of $\z:S\cdot\D@>>>\tO_V$. From (a) we deduce:

(b) {\it The kernel of $\z_0:S\cdot\D^+@>>>\cg$ is equal to $A'$.}
\nl
The image of $\z_0$ is the subgroup of $\cg$ generated by $R$, by $\ty_i\ty_{i'}$ ($i,i'\in I_{odd}$) and by 
$c$. Hence it is equal to $H$. We see that $\z_0$ defines an isomorphism 
$$(S\cdot\D^+)/A'@>\si>>H.\tag c$$
If $I_{odd}\ne\emp$, we can identify 
$$(S\cdot\D^+)/A'=(S^0\cdot\D^+)/A'_1\tag d$$
where 
$$\align&A'_1=\{sc_*^n;s\in S,n\in\ZZ/2,s_i=c'_i{}^{n_i}\text{ if }i\in I_{odd}^{\ge2},\\&s_i=1\text{ if }
i\in I_{even}\cup(I_{odd}-I_{odd}^{\ge2}),n_i\in\ZZ/2,n=\sum_{i\in I_{odd}^{\ge2}}n_i\}.\endalign$$
We show:

(e) {\it If $I_{odd}^{\ge2}=\emp$ then $H=S^0\cdot\D^+$; hence $H/H^0=\D^+$.}
\nl
Assume first that $I_{odd}=\emp$. Then $A'=\{1\}$, $\D^+=\{1,c_*\}$, $S=S^0$ and (c) becomes 
$S^0\sqc S^0c_*@>\si>>H$; thus (e) holds. If $I_{odd}\ne\emp$ and $I_{odd}^{\ge2}=\emp$ then we have 
$A'_1=\{1\}$ so that (e) follows from (d).

\subhead A.7\endsubhead
In this subsection we assume that $I_{odd}^{\ge2}\ne\emp$. 
Using A.6(c),(d) we see that in this case we have $H^0=S^0/A'_2$ where
$$\align&A'_2=\{s\in S;s_i=c'_i{}^{n_i}\text{ if }i\in I_{odd}^{\ge2},s_i=1\\&\text{ if }i\in 
I_{even}\cup(I_{odd}-I_{odd}^{\ge2}),n_i\in\ZZ/2,\sum_{i\in I_{odd}^{\ge2}}n_i=0\},\endalign$$
so that for $i\in I_{odd}^{\ge2}$, the image of $c'_i$ in $H^0$ is independent of $i$; we denote it by $c'$.
Note that the $\D^+$ action on $S^0$ induces an action of $\D^+$ on $S^0/A'_2$ and we can form the 
semidirect product $(S^0/A'_2)\cdot\D^+$. From A.6(c),(d) we see that $H=((S^0/A'_2)\cdot\D^+)/\{1,c'c_*\}$. 
Now $(S^0/A'_2)_{der}$ is the image of $(S^0)_{der}=\prod_{i\in I_{even}\cup I_{odd}^{\ge3}}S_i$ under 
$p:S^0@>>>(S^0/A'_2)_{der}$. The condition that $(S^0/A'_2)_{der}$ is simply connected is that the 
restriction of $p$ to $(S^0)_{der}$ has trivial kernel; that kernel is 
$$\align&A'_2\cap(S^0)_{der}=\{s\in S;s_i=c'_i{}^{n_i}\text{ if }i\in I_{odd}^{\ge3},s_i=1\\&\text{ if }
i\in I_{even}\cup(I_{odd}-I_{odd}^{\ge3}),n_i\in\ZZ/2,\sum_{i\in I_{odd}^{\ge2}}n_i=0\},\endalign$$
and this is trivial precisely when $|I_{odd}^{\ge3}|\le1$.

\subhead A.8\endsubhead
From the results in A.6,A.7 we see that $(H^0)_{der}$ is simply connected if and only if
$|I_{odd}^{\ge3}|\le1$.

\subhead A.9\endsubhead
Let $V$ be a $\CC$-vector space of dimension $N\ge4$ with a fixed symmetric nondegenerate symplectic form 
$(,):V\T V@>>>\CC$. Let $\cg=Is_V$. We assume that $V=\op_{i\ge1}V_i$ where $V_i=W_i\ot E_i$ and $W_i,E_i$ 
are
$\CC$-vector spaces of dimension $i,m_i$ respectively with given nondegenerate bilinear forms $(,)$ (so that
if $i$ is odd, $(,)$ is symmetric for $W_i$ and symplectic for $E_i$; if $i$ is even, $(,)$ is symmetric for
$E_i$ and symplectic for $W_i$) such that $(w\ot e,w'\ot e')=(w,w')(e,e')$ for $w,w'\in W_i$, $e,e'\in E_i$ 
and $(V_i,V_j)=0$ for $i\ne j$. Let $I=\{i\ge1;m_i\ge1\}$, $I_{odd}=I\cap(2\ZZ+1)$, $I_{even}=I\cap(2\ZZ)$. 

We have an imbedding $\t_i:Is_{E_i}@>>>Is_{V_i}$, $g\m 1_{W_i}\ot g$. Let 
$\Ph:\prod_{i\in I}Is_{E_i}@>>>\cg$ be the homomorphism $(g_i)\m\op_{i\in I}\t_i(g_i)$ (an imbedding); it 
restricts to an imbedding $\Ph^0:\prod_{i\in I}Is^0_{E_i}@>>>\cg$. Let $H=\text{ image of }\Ph$. Then $H$ is 
of the form $H(u)$ for some $u\in\cg$ (as in 1.7) and all $H(u)$ as in 1.7 are obtained up to conjugacy. We 
have $\text{ image of }\Ph^0=H^0$. 

We see that $(H^0)_{der}$ is simply connected if and only if the number of $i\in I_{even}$ such that
$m_i\ge3$ is $0$.

\subhead A.10\endsubhead
In this subsection we assume that $L$ is the centralizer of a torus $S$ in $\cg$, that $u\in L$ and that 
$(H^0)_{der}$ is simply connected. Let $\cv'$ be the unipotent radical of $Z_L(u)^0$; let $\cv''$ be the 
unipotent radical of $Z_{L_{der}}(u)^0$. We show:

(a) $Z_L(u)=H'\cv'$ (semidirect product) where $H'$ is reductive and $(H'{}^0)_{der}$ is simply connected;

(b) $Z_{L_{der}}(u)=H''\cv''$ (semidirect product) where $H''$ is reductive and $(H''{}^0)_{der}$ is simply 
connected.
\nl
We can assume that $S\sub H$. We have $Z_L(u)=Z_\cg(u)\cap Z_\cg(S)=HV\cap Z_\cg(S)=Z_H(S)Z_\cv(S)$. Thus we
can assume that $H'=Z_H(S)$ so that $H'{}^0=Z_{H^0}(S)$. Since $(H^0)_{der}$ is simply connected, it follows 
that $(Z_{H^0}(S))_{der}$ is simply connected, proving (a). Let $T$ be the connected centre of $L$. We have 
$L=L_{der}T$. Clearly, $Z_L(u)=Z_{L_{der}}(u)T$. Thus we can assume that $H'=H''T$ so that $H'{}^0,H''{}^0$ 
have the same derived group, proving (b).

\subhead A.11\endsubhead
In the remainder of this appendix we assume that $\cg$ (in 1.7) is of exceptional type and $u,H(u)$ are as in
1.7. As G. Seitz pointed out to the author, the structure of $H=H(u)$ can in principle be extracted from 
\cite{\LS}. More precisely, the tables \cite{\LS, 22.3.1-22.3.5} contain information on the structure of the
Lie algebra of $H^0$, on the structure of $H/H^0$ and the action of $H^0$ on the Lie algebra of $\cg$. From 
this one can recover the precise structure of $H^0$; using in addition the tables \cite{\LS, 22.2.1-22.2.6} 
one can recover the structure of the extension $1@>>>H^0@>>>H@>>>H/H^0@>>>1$. (When $\cg$ is of type $E_7$ 
or $E_6$ then, in addition to the corresponding tables in \cite{\LS} for the corresponding adjoint groups, 
we must use the tables for $E_8$ by regarding $\cg$ as a derived subgroup of a Levi subgroup of a parabolic 
subgroup of a group of type $E_8$; we also use A.10.) In this way we find the following results.

If $\cg$ is of type $E_6$ or $G_2$ then $(H^0)_{der}$ is simply connected.

Assume now that $\cg$ is of type $F_4$. If $u$ is of type $A_1\tA_1$ (notation of \cite{\LS}) then
$H=H^0=PGL_2(\CC)\T SL_2(\CC)$. If $u$ is of type $B_3$ (notation of \cite{\LS}) then $H=H^0=PGL_2(\CC)$.
For all other $u$,  $(H^0)_{der}$ is simply connected.

Assume now that $\cg$ is of type $E_7$. If $u$ is of type $A_2A_1^2$ (notation of \cite{\LS}) then
$H=H^0=SL_2(\CC)^3/\{\pm1\}$ with $\{\pm1\}$ imbedded diagonally in the centre of $SL_2(\CC)^3$.
For all other $u$, $(H^0)_{der}$ is simply connected.

Assume now that $\cg$ is of type $E_8$. If $u$ is of type $A_2A_1^2$ (notation of \cite{\LS}) then
$H=H^0=(SL_2(\CC)\T Spin_7(\CC))/\{\pm1\}$ with $\{\pm1\}$ imbedded diagonally in the centre of
$SL_2(\CC)\T Spin_7(\CC)$. 
If $u$ is of type $A_3A_2A_1$ (notation of \cite{\LS}) then $H=H^0=PGL_2(\CC)\T SL_2(\CC)$.
If $u$ is of type $A_4A_2$ (notation of \cite{\LS}) then $H=H^0=SL_2(\CC)^2/\{\pm1\}$ with $\{\pm1\}$ 
imbedded diagonally in the centre of $SL_2(\CC)^2$.
If $u$ is of type $D_4(a_1)A_2$ (notation of \cite{\LS}) then $H=PGL_3(\CC)\cdot\ZZ/2$ (semidirect product 
with the generator of $\ZZ/2$ acting on $PGL_3(\CC)$ by an outer involution).
If $u$ is of type $D_5(a_1)A_1$ (notation of \cite{\LS}) then $H=H^0=PGL_2(\CC)\T SL_2(\CC)$.
If $u$ is of type $A_6$ (notation of \cite{\LS}) then $H=H^0=SL_2(\CC)^2/\{\pm1\}$ with $\{\pm1\}$ imbedded 
diagonally in the centre of $SL_2(\CC)^2$.
For all other $u$,  $(H^0)_{der}$ is simply connected.

\enddocument